\documentclass[12pt]{article}
\usepackage{amsfonts,amsmath,amssymb,latexsym}
\newtheorem{theorem}{Theorem}

\newcommand{\R}{\mathbb{R}}

\newcommand\bea{\begin{eqnarray*}}
\newcommand\eea{\end{eqnarray*}}
\newcommand\be{\begin{equation}}
\newcommand\ee{\end{equation}}

\def\<{\langle}
\def\>{\rangle}
\newcommand{\po}{{\hspace*{-1ex}}{\bf .  }}
\newcommand\proof{\noindent{\it Proof: }}
\newcommand\qed{\ifhmode\unskip\nobreak\fi\ifmmode\ifinner\else
\hskip5 pt \fi\fi\hbox{\hskip5 pt \vrule width4 pt height6 pt
depth1.5 pt \hskip 1pt }}
\begin{document}

\title{Entire bounded constant mean curvature Killing graphs}
\author{M. Dajczer and J. H. de Lira\thanks
{Partially supported by CNPq and FUNCAP.}}
\date{}
\maketitle

\begin{abstract}
We show that under certain curvature conditions of the ambient space an 
entire Killing graph of constant mean curvature lying inside  a slab                                                             
must be a totally geodesic slice.
\end{abstract}
%\vspace{0.3cm}

%{\small \noindent {\bf Keywords:} Conformal Killing graphs,
%prescribed mean curvature, quasilinear elliptic PDE. \noindent {\bf MSC 2000:}
%53C42, 53A10.}

\section{Introduction}

Let  $N^{n+1}$ denote a complete  Riemannian manifold carrying a non-singular 
Killing vector field $Y$ with integrable orthogonal distribution. Hence 
$N^{n+1}$ is foliated by complete isometric totally geodesic hypersurfaces
(called slices) that we assume to be noncompact. It is a natural problem 
to search for conditions that imply that any entire constant mean curvature graph 
over a slice along the flux of $Y$  must itself be a slice.

In the particular case of a Riemannian product $N^{n+1}=M^n\times\R$ it was 
shown by Rosenberg, Schulze and Spruck \cite{RSS} that if the Ricci 
curvature of $M^n$  satisfies $\textrm{Ric}_M\geq 0$ and its sectional 
curvature satisfies $K_M\geq -K_0$  for a nonnegative constant $K_0$, then any 
entire  minimal graph over $M^n$ with nonnegative height function must be a 
slice.  This result extends the  celebrated  theorem by Bombieri, 
De Georgi and Miranda \cite{BGM} for entire minimal hypersurfaces in Euclidean 
space. In the latter case, it suffices to assume constant mean curvature
due to a result of Chern \cite{Ch} and Flanders \cite{Fl}.

In our more general situation, the ambient space  $N^{n+1}$ has a Riemannian 
warped product structure
$$
N^{n+1}=M^n\times_\rho\R
$$
where the warping function $\rho\in C^\infty(M)$ is $\rho=|Y|>0$.  
\newpage

Let $\Psi\colon\R\times M^n\to N^{n+1}$ denote the flux generated 
by the Killing field $Y$. Fix a integral hypersurface 
$M^n$ of the orthogonal distribution.  Then, the  \emph{Killing graph} $\Sigma(u)$ 
associated to  $u\in C^\infty(M)$ is the  hypersurface defined as
$$
\Sigma(u)=\{\Psi(u(x),x): x\in M^n\}.
$$
That the graph $\Sigma(u)$ is bounded means that it is contained in a slab 
in $N^{n+1}$, that is, there exists a closed interval $\mathbb{I}\subset\R$ 
such that $\Sigma(u)\subset\Psi(\mathbb{I}\times M^n)$. 

\begin{theorem}\label{main}\po Let $N^{n+1}=M^n\times_\rho\R$ be a Riemannian 
warped product where $M^n$ is complete noncompact with $K_M\geq- K_0$ for 
$K_0\geq 0$ and $\textrm{Ric}_M\geq 0$. If
$\rho\geq\rho_0>0$ satisfies  $\|\rho\|_1<\infty$ and if $\textrm{Ric}_N\geq-L$ 
for $L\geq 0$, then any bounded entire  Killing graph with constant mean 
curvature must be a slice.
\end{theorem}

The condition that $\textrm{Ric}_N\geq-L$ for some constant $L\geq 0$ can be 
dropped if we assume  that $\|\rho\|_2<\infty$.  In fact, this follows from 
the relation between the   Ricci curvatures of $N^{n+1}$ and $M^n$ given by 
$$
\textrm{Ric}_N({\sf v}, {\sf v}) =\textrm{Ric}_M(\pi_*{\sf v}, \pi_*{\sf v})
-\frac{1}{\rho}\textrm{Hess}^M\rho(\pi_*{\sf v},\pi_*{\sf v})
-\frac{1}{\rho^3}\<{\sf v},Y\>^2\Delta^M\rho 
$$
where $\pi$ denotes the projection from $N^{n+1}$ to the factor $M^n$.

The key ingredient in the proof of the above result to reduce the case from
constant mean curvature to minimal is the  Omori-Yau maximum principle for 
the Laplacian in the sense of Pigola, Rigoli and Setti given in \cite{PRS}.
For other results on entire constant mean curvature graphs inside slabs 
we refer to \cite{AD}, \cite{RSV} and \cite{Sa}.

\section{A gradient estimate}

In this section, we extend in several directions the gradient estimate for 
minimal graphs in ambient product spaces given by Theorem 4.1 of \cite{RSS}.
For that purpose we make use of the Korevaar-Simon method \cite{Ko} to show 
the existence of an apriori gradient estimate for a nonnegative solution 
of (\ref{pde}) for constant mean curvature $H$ over a geodesic ball $B_R(p)$ 
in $M^n$ centered at $p$ and radius $R$.  \vspace{1ex}

It was shown in \cite{DHL} that a Killing graph $\Sigma(u)$ has mean 
curvature $H$ if and only if $u\in C^2(M)$ satisfies the elliptic PDE of 
divergence form
\be\label{pde}
\textrm{div}_M\bigg(\frac{\nabla^M u}{W}\bigg)
-\frac{1}{2\gamma W}\<\nabla^M\gamma,\nabla^M u\>=nH
\ee
where $\gamma=1/\rho^2$. Here, the mean curvature $H$ is computed with respect 
to the orientation given by the Gauss map
\be\label{gauss}
\mathcal{N}=\frac{1}{W}(\gamma Y-\Psi_*\nabla^M u)\;\;\mbox{where}
\;\;W=\sqrt{\gamma+|\nabla^M u|^2}.
\ee

In the next result the function  $G\in C^\infty([0,+\infty))$  satisfies
$f''=Gf$ where $f\in C^\infty([0,+\infty))$ is positive  outside the origin, 
$f(0)=0$, $f'(0)=1$ and  $g=f'/f$ verifies that $tg(t)$ is bounded in $[0,r]$ for 
any given $r>0$.  For instance, this is the case for $f(t)=\sinh t$.

\begin{theorem}\label{est}\po 
Let $N^{n+1}=M^n\times_\rho\R$ be a complete Riemannian manifold 
satisfying  $\textrm{Ric}_N\geq-L$ for some constant $L\geq 0$.
Assume that the radial sectional curvatures along geodesic issuing 
from a fixed point $p\in M^n$  satisfy  $K_{rad}\geq- G(d)$ where 
$d(x)=\textrm{dist}_M(x,p)$ is the geodesic distance.
Let $\Sigma(u)$ be a Killing graph with  constant mean 
curvature $H$ over $B_R(p)$ with $u\geq 0$. Then,
\be\label{estloc}
|\nabla^M u(p)|\leq  Ce^D
\ee
where the constants $C$ and $D$ are given by (\ref{c0}) and (\ref{c00}).
\end{theorem}

\proof Let $\nabla$ and $\Delta$ always denote  the gradient
and Laplace operator on the graph $\Sigma=\Sigma(u)$, respectively.  For $H$ constant
it is  well-known \cite{FR} that
$$
\Delta\< Y,\mathcal{N}\> 
= -(|A|^2+\textrm{Ric}_N(\mathcal{N},\mathcal{N}))\<Y,\mathcal{N}\>
$$
where $|A|$ denotes de norm of the second fundamental form of the graph.
Since $\<Y,\mathcal{N}\> = 1/W$, we obtain 
\be\label{four}
\Delta W-\frac{2}{W}|\nabla W|^2=(|A|^2+\textrm{Ric}_N(\mathcal{N},\mathcal{N}))W.
\ee

Let $\mathcal{K}^+$ denote the solid half-cylinder $\Psi(\R^+\times B_R)$
where $B_R=B_R(p)$. 
The Killing field is $Y=\partial/\partial s$ where $s$ parametrizes 
the line factor in $N^{n+1}$. Define a function $\phi:\mathcal{K}^+\to \R$ by
$$
\phi(y)=\left(1-\frac{d^2(x)}{R^2}-C_0s\right)^+
\quad \textrm{if} \quad  y=\Psi(s,x) \in K^+
$$
where $C_0=1/2u(p)>0$ and $+$ means the positive part.
Now consider the function $\eta$ supported in $B_R$ given by
\be\label{def}
\eta(y) = e^{K\phi(y)}-1, \quad y=\Psi(u(x),x),
\ee
for same constant $K>0$.

Set $h=\eta W$ and let $q\in B_R$ be a (necessarily interior) point of maximum of $h$.
First consider the case when $q \in B_R\backslash C(p)$, where $C(p)$ 
is the cut locus of $p$ in $M^n$, so $h$ is smooth near $q$.  
Without further reference we compute at the point $\Psi(u(q),q)$.  
It holds that
$$
\nabla h=\eta\nabla W+W\nabla\eta=0
$$
and
$$
\Delta h=W\Delta\eta+\Big(\Delta W-\frac{2}{W}|\nabla W|^2\Big)\eta.
$$
Since the Hessian form of $h$ is nonpositive, 
we have from (\ref{four}) that
$$
\Delta h=W\big(\Delta\eta +(|A|^2
+\textrm{Ric}_N(\mathcal{N},\mathcal{N}))\eta\big)\leq 0.
$$
The Ricci curvature assumption now yields 
$\Delta\eta\leq L\eta$
and from (\ref{def}) we have 
\be\label{zz}
\Delta\phi + K|\nabla \phi|^2\leq \frac{L}{K}.
\ee

In the sequel, we estimate both terms in the left hand side of (\ref{zz}). 
We have that $\bar\nabla s=\gamma Y$  where $\bar\nabla$ denote the Riemannian 
connection in the ambient space.
Where $\phi$ is differentiable and positive, we have
\be\label{gradphi}
\bar\nabla \phi (y)=-\frac{2d}{R^2}\bar\nabla d (y)-C_0\gamma Y(y)
=-\frac{2d}{R^2}\Psi_*(s,x)\nabla^{M} d(x)-C_0\gamma Y(y).
\ee
Then,
\be\label{old}
\<\bar\nabla\phi,\mathcal{N}\>=\frac{2d}{R^2W}\<\nabla^Md,\nabla^Mu\>
-C_0\frac{\gamma}{W},
\ee
and thus
\bea
|\nabla\phi|^2\!\!\!&=&\!\!\! |\bar\nabla \phi|^2-\<\bar\nabla\phi,\mathcal{N}\>^2\\
\!\!\!&=&\!\!\!\frac{4d^2}{R^4}\Big(1-\frac{1}{W^2}\<\nabla^Md,\nabla^Mu\>^2\Big) 
+\frac{C_0^2\gamma}{W^2}|\nabla^Mu|^2
+\frac{4dC_0\gamma}{R^2W^2}\<\nabla^Md,\nabla^Mu\>\\
\!\!\!&\geq &\!\!\!\frac{4d^2}{R^4}\Big(1-\frac{1}{W^2}|\nabla^M u|^2\Big)
+\frac{C_0^2\gamma}{W^2}|\nabla^Mu|^2-\frac{4dC_0\gamma}{R^2W^2}|\nabla^M u|\\
\!\!\!&\geq&\!\!\!  \gamma\Big(\frac{C_0}{W}|\nabla^M u|
-\frac{2}{RW}\Big)^2
\eea
where to obtain the last inequality we used that $d\leq R$. 

Observe that having
\be\label{ine}
\left(\frac{C_0}{W}|\nabla^M u|-\frac{2}{RW}\right)^2\geq\frac{C_0^2}{\alpha^2}
\ee
for some constant $\alpha>2$ is equivalent to
$$
\left(\frac{C_0}{W}|\nabla^M u|-\frac{2}{RW}
-\frac{C_0}{\alpha}\right)\left(\frac{C_0}{W}|\nabla^M u|
-\frac{2}{RW}+\frac{C_0}{\alpha}\right) \ge 0.
$$
Clearly,  the latter holds if the first factor is nonnegative
or, equivalently, if
$$
\alpha|\nabla^M u|-W\ge\frac{2\alpha}{C_0R}.
$$
On the other hand,  that
$$
\alpha|\nabla^M u|-W\ge W
$$
is equivalent to
$$
\alpha^2|\nabla^M u|^2\ge 4W^2\;\;\mbox{or}\;\;
|\nabla^M u|^2\ge 4\gamma/(\alpha^2-4)\;\;\mbox{or}\;\;
W^2\ge\frac{\alpha^2}{(\alpha^2-4)\rho^2}.
$$
Therefore,  we assume that
\be\label{est-contrad}
W(q)\geq \max\left\{\frac{\alpha}{\min_{B_R}\rho\sqrt{\alpha^2-4}}, 
\frac{2\alpha}{C_0R}\right\}
\ee
holds and obtain that (\ref{ine}) also holds.  Hence, assuming (\ref{est-contrad}) 
we obtain that
\be\label{first}
|\nabla\phi|^2\geq \frac{C_0^2}{\alpha^2\rho^2}=\frac{1}{4\alpha^2\rho^2u^2(p)}.
\ee

Now we estimate the other term in (\ref{zz}). If  $\{e_i\}_{i=1}^n$ is a local 
orthonormal frame in $\Sigma$, we have from
$$
\Delta \phi=\sum_{i=1}^n\< \nabla_{e_i}\nabla\phi, e_i\> 
= \sum_{i=1}^n\< \bar\nabla_{e_i}(\bar\nabla\phi
-\<\bar\nabla\phi,\mathcal{N}\>\mathcal{N}), e_i\>
$$
that
\be\label{one}
\Delta \phi=\sum_{i=1}^n \<\bar\nabla_{e_i}\bar\nabla \phi,e_i\>
+nH\<\bar\nabla\phi,\mathcal{N}\>.
\ee
We obtain from (\ref{gradphi}) that
\be\label{new}
\<\bar\nabla_{e_i}\bar\nabla\phi,e_i\>=-\frac{2}{R^2}(d\<\bar\nabla_{e_i}
\bar\nabla d,e_i\>+\<\bar\nabla d,e_i\>^2)-C_0\<\bar\nabla_{e_i}\gamma Y,e_i\>.
\ee
On the other hand,
$$
\<\bar\nabla_{e_i}\bar\nabla d,e_i\>=\textrm{Hess}^M d(\hat{e}_i,\hat{e}_i)
+\gamma^2\<Y,e_i\>^2\<\bar\nabla_Y\bar\nabla d,Y\>
$$
where $\hat{e}_i = \pi_* e_i$ and $\pi$ is the projection from $N^{n+1}$ to $M^n$.
Our assumptions and the comparison Theorem 2.3  for the Hessian in \cite{PRS2} yield 
$$
\textrm{Hess}^M d \leq g(d)(\<\cdot,\cdot\>-\textrm{d}d\otimes\textrm{d}d)
$$
where $g=f'/f$. Therefore, using that for the warped product metric
$$
\bar\nabla_Y\bar\nabla d= \frac{1}{\rho}\<\nabla^Md,\nabla^M\rho\>Y,
$$
we obtain that
\be\label{laplace}
\<\bar\nabla_{e_i}\bar\nabla d,e_i\>\leq 
g(\<\hat e_i,\hat e_i\>-\<\nabla^M d,\hat e_i\>^2)
+\frac{1}{2}\<Y,e_i\>^2|\nabla^M\gamma|.
\ee
It follows from (\ref{new}) that
\bea
\<\bar\nabla_{e_i}\bar\nabla \phi, e_i\> 
\!\!\!&\geq&\!\!\! -\frac{2d}{R^2}(g(1-\gamma\<Y,e_i\>^2-\<\nabla^M d,e_i\>^2)
+\frac{1}{2}\<Y,e_i\>^2|\nabla^M\gamma|)\\
\!\!\!&&\!\!\!-\frac{2}{R^2}\<\bar\nabla d, e_i\>^2-C_0\<\bar\nabla_{e_i}\gamma Y,e_i\>.
\eea
Since $\<\bar\nabla_Y Y,Y\>=0$ and hence $Y(\gamma)=0$, we have
\bea
\<\bar\nabla_{e_i}\gamma Y,e_i\> \!\!\!&=&\!\!\!\<\bar\nabla_{\hat{e}_i}\gamma Y,\hat{e}_i\>
+ \gamma \< e_i,Y\>\<\bar\nabla_Y\gamma Y,\hat{e}_i\>
+ \gamma\<e_i,Y\>\<\bar\nabla_{\hat{e}_i}\gamma Y,Y\>\\
\!\!\!&=&\!\!\!\gamma^2\<e_i,Y\>\<\bar\nabla_YY,\hat{e}_i\>
+\gamma\<e_i,Y\>(\hat{e}_i(\gamma)\<Y,Y\>
+\gamma\<\bar\nabla_{\hat{e}_i}Y,Y\>)\\
\!\!\!&=&\!\!\!\gamma^2\<e_i,Y\>\<\bar\nabla_YY,\hat{e}_i\>
+\gamma\<e_i,Y\>(\hat{e}_i(\gamma)\<Y,Y\>
-\gamma\<\bar\nabla_YY,\hat{e}_i\>)\\
\!\!\!&=&\!\!\!\<e_i,Y\>e_i(\gamma).
\eea
We obtain that
\bea
\sum_{i=1}^n \<\bar\nabla_{e_i}\bar\nabla \phi, e_i\>\geq \!\!\!&-&\!\!\!
\frac{2d}{R^2}\Big(g(n-\gamma|Y^T|^2- |\nabla^\Sigma d|^2)
+ \frac{1}{2}|\nabla^M\gamma||Y^T|^2\Big)\\
\!\!\!&-&\!\!\!\frac{2}{R^2}|\nabla^\Sigma d|^2-C_0\<\bar\nabla\gamma, Y^T\>
\eea
where $Y^T$ denotes the component of $Y$ tangent to $\Sigma$. Since
$$
\<\bar\nabla\gamma, Y^T\>
=-\<Y,\mathcal{N}\>\<\bar\nabla\gamma, \mathcal{N}\>
= -\frac{1}{W}\<\bar\nabla\gamma,\mathcal{N}\>,
$$
it follows that 
$$
\Delta\phi\ge-\frac{2d}{R^2}(ng+(\kappa-g)\gamma|Y^T|^2)
+\frac{2}{R^2}(dg-1)|\nabla^\Sigma d|^2 
+\frac{C_0}{W}\<\bar\nabla\gamma,\mathcal{N}\>+nH\<\bar\nabla\phi,\mathcal{N}\>
$$
where 
$$
\kappa=\frac{1}{\rho}|\nabla^M\rho|=\frac{1}{2\gamma}|\nabla^M\gamma|.
$$
We have from (\ref{old}) that
$$
nH\<\bar\nabla\phi,\mathcal{N}\>=
-\frac{2d}{R^2}\<\bar\nabla d,\mathcal{N}\>nH-\frac{C_0\gamma}{W}nH.
$$
Moreover, 
$$
\frac{C_0}{W}(\<\bar\nabla\gamma,\mathcal{N}\>-n\gamma H)
\geq -\frac{\gamma C_0}{W}(2\kappa+n|H|)\ge - C_0\sqrt{\gamma}(2\kappa+n|H|).
$$
We obtain that
$$
\Delta\phi\geq-\frac{2d}{R^2}(ng+(\kappa-g)\gamma|Y^T|^2+n|H|)
+\frac{2}{R^2}(dg-1)|\nabla^\Sigma d|^2
-C_0\sqrt{\gamma}(2\kappa+n|H|).
$$
Thus,
$$
\Delta\phi\geq-\frac{2}{R^2}(1+ndg+\kappa+n|H|)
-\frac{1}{2\rho u(p)}(2\kappa+n|H|).
$$
We denote $g_0=\max_{[0,R]}dg$, $\kappa_0=\max_{B_R}\kappa$, $\rho_1
=\min_{B_R}\rho$ and 
\be\label{ab}
A_0=1+ng_0+\kappa_0+n|H|,\;\;\;A_1=\frac{1}{4\rho_1}(2\kappa_0+n|H|).
\ee
It follows that
\be\label{est2}
\Delta\phi>-2E\;\;\;\mbox{where}\;\;\;E=\frac{A_0}{R^2}+\frac{A_1}{u(p)}.
\ee
Using  (\ref{first}) and (\ref{est2}) and choosing $\alpha=4$, 
we obtain from (\ref{zz}) that
\be\label{zzz}
\frac{L}{K}\geq \Delta\phi + K|\nabla \phi|^2>-2E+KP.
\ee
where $P=1/64\delta^2u^2(p)$ and  $\delta=\max_{B_R}\rho$.
Taking
\be\label{k0}
K>\frac{1}{P}\left(E+\sqrt{E^2+PL}\right)
\ee
we obtain a contradiction in (\ref{zzz}) and  
conclude that (\ref{est-contrad}) cannot hold. Thus,
\be\label{c0}
W(q)<C=\max\left\{\frac{2}{\sqrt{3}\,\rho_1},
\frac{16u(p)}{R}\right\}.
\ee
Hence,
$$
h(p)=(e^{K/2}-1)W(p)\leq h(q)\leq C(e^K-1)
$$
where $K$ is given by (\ref{k0}). We obtain (\ref{estloc}) choosing 
\be\label{c00}
D=\log (e^{K/2}+1).
\ee

To conclude the proof we observe that the same argument in \cite{RSS} 
proves that the restriction that  $q\not\in C(p)$ can be dropped. \qed

\section{The proof of Theorem \ref{main}}

First we obtain a gradient estimate for the entire graph $\Sigma(u)$
taking $R\to +\infty$ in Theorem \ref{est}. Notice that we need 
$g_0=\sup_{[0,+\infty]}tg(t)<+\infty$
to hold. But this is the case choosing
$$
f(t)=\frac{1}{\sqrt{K_0}}\sinh\sqrt{K_0}\,t.
$$
In fact, the functions $f$, $g=\sqrt{K_0}\coth\sqrt{K_0}\,t$ and $G=K_0$ 
now meet all of the requirements in Theorem \ref{est}. 

Taking $R\to +\infty$ we conclude that we can choose
$$
C=\frac{2}{\sqrt{3}\,\delta}\;\;\mbox{and}\;\;
E=\frac{1}{4\delta u(p)}(2\kappa_0+n|H|)
$$
where now $\kappa_0=\max_M\kappa$. Then, we can take
$$
K=8\delta(A+\sqrt{A^2 + L}) u(p)\;\;\mbox{where}\;\; A=2(2\kappa_0+n|H|).
$$
Up to a vertical translation, we may assume 
that $\mathbb{I}=[0,T]$ for some constant $T>0$. Thus, we can take 
$$
K=8\delta(A+\sqrt{A^2 + L})T
$$ 
and conclude that there exists a uniform bound for the gradient, i.e.,
\be\label{bound}
|\nabla^M u|<C_1
\ee
at any point for some constant $C_1>0$ independent of $u$. 

Next, we argue that $H=0$.
We claim that the Omori-Yau maximum principle for the Laplacian holds 
on $\Sigma(u)$ in the sense of Pigola, Rigoli and Setti in \cite{PRS}.  
For that we show that the nonnegative function $\psi\in C^\infty(\Sigma)$ 
defined by $\psi(x,u(x))=d^2(x)$  satisfies the three conditions that are 
required  by Theorem 1.9 in \cite{PRS}.  The first holds since 
$\psi(x)\to +\infty$ as $x\to\infty$. Also the second is satisfied since
$$
|\nabla\psi|=2d|\nabla d|\leq 2d|\nabla^M d|=2\sqrt{\psi}.
$$
For the final condition,  we have 
$$
\Delta\psi=2d\Delta d + 2|\nabla d|^2\leq 2(d\Delta d + 1).
$$
We obtain from (\ref{one}) that
$$
\Delta d\leq\sum_{i=1}^n \<\bar\nabla_{e_i}\bar\nabla d,e_i\>+n|H|
$$
and from (\ref{laplace}) that
$$
\<\bar\nabla_{e_i}\bar\nabla d,e_i\>\leq \coth\sqrt{K_0}\,d
+\frac{1}{\rho}|\nabla^M\rho|.
$$
We conclude that $\Delta\psi\leq C_2\sqrt{\psi}$ for some positive 
constant $C_2$ outside a compact set.  This implies that the third condition
is also satisfied and concludes the proof of the claim.

Since $u=s|_\Sigma$ and $\bar\nabla s = \gamma Y$,  we have 
\bea
\Delta u \!\!\!&=&\!\!\! \sum_{i=1}^n\<\nabla_{e_i}\nabla s,e_i\> 
=\sum_{i=1}^n\<\nabla_{e_i}\gamma Y^T,e_i\>\\
\!\!\!&=&\!\!\!\gamma\sum_{i=1}^n\<\bar\nabla_{e_i}(Y-\<Y,\mathcal{N}\>\mathcal{N}),e_i\>
+\frac{1}{\gamma}\sum_{i=1}^n\<e_i,\nabla\gamma\>\<e_i,\bar\nabla s\>\\
\!\!\!&=&\!\!\!n\gamma H\<Y,\mathcal{N}\>+\frac{1}{\gamma}\<\nabla\gamma,\nabla u\>.
\eea
Thus,
\be\label{final}
\Delta u =\frac{n\gamma H}{\sqrt{\gamma+|\nabla^M u|^2}}
+\frac{1}{\gamma}\<\nabla\gamma,\nabla u\>.
\ee

We apply the maximum principle for the Laplacian to $u$ as a bounded 
function on $\Sigma$.  Since $u$ is bounded from above there exists 
a sequence $\{x_k\}_{k\in\mathbb{N}}$ such that
\be\label{first2}
|\nabla u(x_k)|< 1/k\;\;\; \mbox{and}\;\;\;\;\Delta u(x_k)<1/k.
\ee
Since $u$ is bounded from below there is a sequence 
$\{y_k\}_{k\in\mathbb{N}}$ such that
\be\label{second}
|\nabla u(y_k)|< 1/k\;\;\; \mbox{and}\;\;\;\;\Delta u(y_k)>-1/k.
\ee
It follows from (\ref{bound}),  (\ref{final}), (\ref{first2}) and 
(\ref{second}) that $H=0$.

To conclude the proof, we observe that now (\ref{final}) reduces to
$$
\Delta u +\frac{2}{\rho}\<\nabla\rho,\nabla u\>=0.
$$
This can be written as $Lu=0$ where
$$
Lu=e^{-\varphi}\textrm{div}_\Sigma e^{\varphi}\nabla u
$$
and $\varphi=2\log\rho$.
Let $\{x_i\}_{i=1}^n$ be a local system of coordinates in $M^n$ with
metric $\sigma_{ij}=\<\partial x_i,\partial x_j\>$.  Then, we have
$$
L=e^{-\varphi}\textrm{div}_\Sigma (e^\varphi g^{ij}u_i\partial x_j)
$$
where 
$$
g^{ij}=\sigma^{ij}-\frac{u^iu^j}{W^2}\;\;\mbox{and}\;\;u^i=\sigma^{ik}u_k.
$$
It follows from (\ref{bound}) that $L$ is an uniformly elliptic operator
in divergence form.  Now, if we view $L$ as a operator acting on $M^n$ 
and since $\textrm{Ric}_M\geq 0$, it follows from 
Theorem 7.4 in \cite{SC} that $u$ must be constant.\qed

{\renewcommand{\baselinestretch}{1} \hspace*{-20ex}\begin{tabbing}
\indent \= Marcos Dajczer\\
\> IMPA \\
\> Estrada Dona Castorina, 110\\
\> 22460-320 -- Rio de Janeiro -- Brazil\\
\> marcos@impa.br\\
\end{tabbing}}

\vspace*{-4ex}

{\renewcommand{\baselinestretch}{1} \hspace*{-20ex}\begin{tabbing}
\indent \= Jorge Herbert S. de Lira\\
\> Departamento de Matematica - UFC, \\
\> Bloco 914 -- Campus do Pici\\
\> 60455-760 -- Fortaleza -- Ceara -- Brazil\\
\> jorge.lira@pq.cnpq.br
\end{tabbing}}

\end{document}